\documentclass{amsart} 

\usepackage{amsfonts}
\usepackage{amsmath}
\usepackage{amsthm}

\newcommand{\Q}{\mathbb{Q}_p}
\newcommand{\Z}{\mathbb{Z}}
\newcommand{\la}{\lambda}

\newcommand{\ar}{\right|}
\newcommand{\al}{\left|}
\newtheorem{vd1}{Lemma}
\newtheorem{vd}[vd1]{Lemma}
\newtheorem{cor}[vd1]{Corollary}
\theoremstyle{plain}
\newtheorem{res}{Lemma}
\newtheorem{res2}[res]{Corollary}
\newtheorem{vdc}[res]{Corollary}
\newtheorem{one}[res]{Lemma}
\newtheorem{two}[res]{Lemma}

\newtheorem{four}[res]{Lemma}

\newtheorem{threetwo}[res]{Lemma}

\usepackage{amscd}

\begin{document}

\address{School of Mathematics,
University of New South Wales,
Sydney, NSW 2052,
Australia}
\email{keith@maths.unsw.edu.au}
\subjclass[2000]{Primary 43A70, Secondary 11F85}
\author{Keith Rogers}
\date{}
\title{A van der Corput Lemma for the $p$-Adic Numbers}
\keywords{Van der Corput lemma, $p$-adic, oscillatory integrals}

\begin{abstract}We prove a version of van der Corput's Lemma for polynomials over the $p$-adic numbers. 
\end{abstract}

\maketitle

\section{Introduction}\label{vdcls}

The following lemma goes back to J.G. van der Corput in \cite{van}. It has many applications in number theory and harmonic analysis. In particular, it is key to
the study of oscillatory integrals (see \cite{st}). We note that only partial van der
Corput type lemmas are known in dimensions greater than one (see \cite{ca}). As a consequence the theory of
oscillatory integrals in higher dimensions is relatively open. 
\begin{vd}\label{vd} Suppose that $f:(a,b)\to \mathbb{R}$ is $n$ times differentiable, where $n\ge2,$ and $|f^{(n)}(x)|\ge \la >0$
on $(a,b)$. Then  
$$
\left|\int_a^be^{ if(x)}\,dx\right| \le C_n\frac{n}{\la^{1/n}},
$$
where $C_n\le 2^{5/3}$ for all $n\ge 2$ and $C_n\to 4/e$ as $n\to \infty.$
\end{vd}

It can be shown by considering $f(x)=x^n$ that the linear growth in $n$ is
optimal, and this more precise formulation is due to G.I. Arhipov, A.A.
Karacuba and V.N. Cubarikov 
\cite{ar}.  Consideration of the Chebyshev polynomials shows that the constant becomes sharp as $n$ tends to infinity (see \cite{me}). 
The following corollary can be easily obtained using Stirling's formula.
\begin{cor}\label{cor} Let $f(x)=a_0+a_1x+\dots+a_nx^n$ be a real polynomial of degree $n\ge 1.$ Then 
$$
\al \int_a^be^{if(x)}\,dx\ar \le \frac{2^{5/3}e}{|a_n|^{1/n}} < \frac{9}{|a_n|^{1/n}}
$$
for all $a,b\in \mathbb{R}$.
\end{cor}
We will prove a $p$-adic version of this corollary, opening the way for
the study of oscillatory integrals on the $p$-adics.
This problem was first considered by J.
Wright \cite{man}, where lemmas for polynomials of degree two and monomials of degree three
were proven.

\section{Introduction to the $p$-Adic Numbers} For a more complete
introduction to the $p$-adic numbers, see \cite{ge} or \cite{sa}. Here we will outline what we
will need. 

Fix a prime number $p$. Any non-zero rational number $x$ can be uniquely expressed in the form $p^{k}m/n,$ where $m$ and $n$ have no common divisors and neither is
divisible by $p$. We then define the $p$-adic norm on the rational numbers by $|x|=p^{-k}$ when $x\neq 0$, and $|0|=0$. We obtain the
$p$-adic numbers by completing $\mathbb{Q}$ with respect to this norm. 
It is not difficult to show that the norm
satisfies the following properties:
\begin{eqnarray*}
& |xy|=|x||y|& \\
&|x+y|\le \max\{|x|,|y|\}.&
\end{eqnarray*}
It follows from the second property that
$$
\{y:|y-x_1|\le p^r\}=\{y:|y-x_2|\le p^r\}$$
when $|x_1-x_2|\le p^r$, so every point within a ball can be considered to be its centre.

A nonzero $p$-adic number $x$ with $|x|=p^{-k},$ may be written in the form
$$
x=\sum_{j=k}^{\infty}x_jp^j,
$$
where $0\le x_j\le p-1$, and $x_k\neq 0.$ This will be called the standard $p$-adic expansion and the arithmetic is done
formally with carrying.   
Define $\chi : \Q \to \mathbb{C}$ by
$$
\chi(x) = \left\{ \begin{array}{ll}
\prod_{j=k}^{-1}e^{2\pi i x_j/p^j} & |x|>1\\
1 & |x|\le 1.
\end{array} \right.
$$
The characters of $\Q$ are all of the form $\chi_{\epsilon}:\Q\to \mathbb{C}$;
$$
\chi_{\epsilon}(x)=\chi(\epsilon x), 
$$
where $\epsilon\in\Q.$ 
Finally $\Q$ is a locally compact commutative group, so there is a Haar measure, that necessarily satisfies $d(ax)=|a|d
x.$ We normalise the measure so that $\{x\in\Q : |x|\le p^r\}$ has measure $p^r.$ 

The usual arguments can be employed to obtain the standard Fourier results.

\section{$p$-Adic van der Corput Lemmas}\label{weeman}

The main thrust of this section is to prove the following lemma for $p$-adic polynomials. The Euclidean arguments will not be helpful as there is no order on the $p$-adic numbers.
\begin{res}\label{res}Suppose that $a_0, \ldots, a_n \in \Q.$ Then 
$$
\al\int_{|x|\le 1} \chi(a_1x+\dots+ a_nx^n)\,dx \ar \le \frac{p^{m}}{|ma_m|^{1/m}},
$$ 
where $m=\max\{l : |la_l|\ge|ja_j| \ \textrm{ for all } \ j\neq l\}.$
\end{res}
Before proving Lemma~\ref{res} we note some easy corollaries.
\begin{res2}\label{res2}Suppose that $a_0, \ldots, a_n \in \Q.$ Then 
$$
\al\int_{|x|\le 1} \chi(a_1x+\dots+ a_nx^n)\,dx \ar \le \frac{2p^{n}}{\lambda^{1/n}},
$$ 
where $\lambda=\max_{1\le j \le n}|a_j|.$
\end{res2}
\begin{proof} Suppose that $|a_k|=\max_{1\le j \le n}|a_j|=\lambda.$ By Lemma~\ref{res} we have
$$
|I|= \al\int_{|x|\le 1} \chi(a_1x+\dots+ a_nx^n)\,dx \ar \le\frac{p^{m}}{|ma_m|^{1/m}},
$$
where $m=\max\{l : |la_l|\ge|ja_j| \
\textrm{ for all } \ j\neq l\}.$ 
Now as $|ma_m|\ge |ka_k|,$ we have
$$
|I| \le\frac{p^{m}}{|ka_k|^{1/m}}\le\frac{k^{1/m}p^{m}}{|a_k|^{1/m}} \le
\frac{n^{1/m}p^{m}}{\lambda^{1/n}}\le \frac{2p^n}{\lambda^{1/n}},
$$
and we are done.
\end{proof}

From this we obtain our main result which holds uniformly over all balls. It is the $p$-adic equivalent of
Corollary~\ref{cor} in Section~\ref{vdcls}. 
\begin{vdc}\label{vdc} Suppose that $x_0,a_0, \ldots, a_n \in \Q$ and $r\in \Z.$ Then 
$$
\al\int_{|x-x_0|\le p^r} \chi(a_0+a_1x+\dots+ a_nx^n)\,dx \ar \le \frac{2p^{n}}{|a_n|^{1/n}}.
$$ 
\end{vdc}
\begin{proof}
Let $y=p^r(x-x_0),$ so that
\begin{align*}
I&=\int_{|x-x_0|\le p^r} \chi(a_0+a_1x+\dots+ a_nx^n)\,dx \\
&=\int_{|y|\le 1} \chi\left(a_0+a_1\left(\frac{y}{p^r}+x_0\right)+\dots+
a_n\left(\frac{y}{p^{r}}+x_0\right)^n\right)\,dx \\
&=p^r\int_{|y|\le 1} \chi\left(b_0(x_0)+\dots+ \frac{b_{n-1}(x_0)y^{n-1}}{p^{(n-1)r}}+
\frac{a_ny^n}{p^{nr}}\right)\,dy\\
&=:p^rI_1,
\end{align*}
where $b_j(x_0)=a_j + {j+1\choose j}a_{j+1}x_0+\dots+{n\choose j}a_nx_0^{n-j}.$
We also note that
\begin{align*}
|I_1|&=\left|\chi(b_0(x_0))\int_{|y|\le 1} \chi\left(\frac{b_1(x_0)y}{p^r}+\dots+ \frac{b_{n-1}(x_0)y^{n-1}}{p^{(n-1)r}}+
\frac{a_ny^n}{p^{nr}}\right)\,dy\right|\\
&=\left|\int_{|y|\le 1} \chi\left(\frac{b_1(x_0)y}{p^r}+\dots+ \frac{b_{n-1}(x_0)y^{n-1}}{p^{(n-1)r}}+
\frac{a_ny^n}{p^{nr}}\right)\,dy\right|\\
&=:|I_2|.
\end{align*}
Thus
$$
|I|=p^r|I_2|\le p^r\frac{2p^{n}}{|p^{-nr}a_n|^{1/n}}=\frac{2p^{n}}{|a_n|^{1/n}},
$$
by Corollary~\ref{res2}. 
\end{proof}

We now turn to the proof of Lemma~\ref{res}. We will need some preliminary lemmas. Our starting point is a consequence of the fact that balls in $\Q$ have multiple centres.
\begin{one}\label{one} Suppose that $a\in\Q$ and $|a|>1.$ Then
$$
\int_{|x|\le 1} \chi(ax)\,dx =0.
$$
\end{one}

\begin{proof}
First we consider the standard expansion of $a,$ so that
$$a=\sum^{\infty}_{j=-k}a_jp^j,$$ where $k\ge 1$ and $a_{-k}\neq 0$. Now as $$\{x: |x|\le 1\} = \{ x :
|x-p^{k-1}|\le 1\},$$ we have
\begin{equation*}
I=\int_{|x|\le 1} \chi(ax)\,dx  =\int_{|x-p^{k-1}|\le 1} \chi(ax)\,dx.\end{equation*}
If we let $y=x-p^{k-1},$ we see that
$$
I= \int_{|y|\le 1} \chi(a(y+p^{k-1}))\,dy=\chi(ap^{k-1})\int_{|y|\le 1} \chi(ay)\,dy,$$ 
so that
$$I=\chi(ap^{k-1})I.$$
Now as $\chi(ap^{k-1})=e^{2\pi i a_{-k}/p}\neq 1,$ we see that $I=0$.
\end{proof} 

If we let $f(y)=a_0+a_1y+\dots+a_ny^n,$ then we denote
\begin{equation}\label{star}
 b_j(y)=\frac{f^j(y)}{j!}=a_j + {j+1\choose j}a_{j+1}y+\dots+{n\choose j}a_ny^{n-j}.\end{equation}
We will use this notation throughout.
\begin{threetwo}\label{threetwo}
Suppose that $|ma_m|>|ja_j|$ for all $j>m,$ and $|y|\le 1$.  Then 
$$
|mb_m(y)|=|ma_m|>|jb_j(y)|
$$
for all $j>m,$ where $b_j$ is given by $(\ref{star})$. 
\end{threetwo}
\begin{proof}
Suppose that $|ma_m|>|ja_j|$ for all $j>m.$ Then
$$
|ma_m|>\left|{j-1\choose m-1}\right||ja_j|,
$$
so that
$$
|a_m|>\left|{j\choose m}a_j\right|
$$
for all $j>m.$ Thus
$$
|mb_m(y)|=|m|\left|a_m + {m+1\choose m}a_{m+1}y+\dots+{n\choose m}a_ny^{n-m}\right|=|ma_m|
$$
for all $|y|\le1$. Similarly, if $k>j>m$, then
$$
|ma_m|> \left|{k-1\choose j-1}\right||ka_k|,
$$
so that
$$
|ma_m|>\left|j{k\choose j}a_k\right|.$$
Putting these together, 
$$
|mb_m(y)|=|ma_m|> \left|ja_j + j{j+1\choose j}a_{j+1}y+\dots+j{n\choose j}a_ny^{n-j}\right|=|jb_j(y)|
$$
for all $|y|\le 1$.
\end{proof}

\begin{two}\label{two} Suppose that $|a_1|> p$ and $|a_1|>|ja_j|$ for $j>1.$ Then
$$
\int_{|x|\le 1} \chi(a_1x+\dots+ a_nx^n)\,dx = 0.
$$

\end{two}

\begin{proof}
Let $|a_1|=p^{k+1}$ where $k\ge 1.$ We split the integral into $p^{k}$ pieces, so that 
$$
I= \sum_{y=0}^{p^k-1}\int_{|h|\leq p^{-k}}\chi(a_1(y+h)+\dots+a_n(y+h)^n)\,dh.
$$
Now 
$$
I=\sum_{y=0}^{p^k-1}\chi(a_1y+\dots+ a_ny^n)I_1(y),
$$
where 
\begin{align*}
I_1(y)&= \int_{|h|\le p^{-k}}\chi(b_1(y)h+\dots+b_{n-1}(y)h^{n-1}+a_nh^n)\,dh\\
&= \frac{1}{p^k}\int_{|x|\le
1}\chi(b_1(y)p^{k}x+\dots+b_{n-1}(y)p^{(n-1)k}x^{n-1}+a_np^{nk}x^n)\,dx,
\end{align*}
and $b_j$ is given by $(\ref{star})$. 
When $|y|\le 1$, we have
$$
|b_1(y)|=|a_1|>|jb_j(y)|
$$
for all $j>1,$ by Lemma~\ref{threetwo}. Hence
$$
|b_1(y)p^k|=\frac{|a_1|}{p^k}=p,
$$
and
$$
|jb_j(y)p^{jk}|p^{jk}< |b_1(y)p^k|p^k=p^{k+1}.
$$ So if $j > 1,$ then
$$
|b_j(y)p^{jk}|\le \frac{1}{|j|p^{(j-1)k}}\le\frac{j}{p^{(j-1)k}}\le\frac{2}{2^{(2-1)1}}=1.
$$
Thus by Lemma~\ref{one},
\begin{align*}
I_1(y)&=\frac{1}{p^k}\int_{|x|\le
1}\chi(b_1(y)p^{k}x)\chi(b_2(y)p^{2k}+\dots+a_np^{nk}x^n)\,dx\\
&=\frac{1}{p^k}\int_{|x|\le
1}\chi(b_1(y)p^{k}x)\,dx=0
\end{align*}
for all $|y|\le1,$ and we are done. 
\end{proof} 

\begin{four}\label{four} Suppose 
that $|ma_m|> p^2$ and $|ma_m|>
|ja_j|$ for all $j\neq m.$
 Then
$$ 
\int_{|x|\le 1} \chi(a_1x+\dots+ a_nx^n)\,dx = \frac{1}{p}\int_{|x|\le 1} \chi(a_1px+\dots+a_np^nx^n)\,dx.
$$
\end{four}

\begin{proof}
We split the integral into $p$ pieces, so that
$$
I=\int_{|x|\le 1} \chi(a_1x+\dots+ a_nx^n)\,dx =\sum^{p-1}_{y=0}\chi(a_1y+\dots+a_ny^n)I_1(y),
$$
where 
\begin{align*}
I_1(y)&= \int_{|h|\le 1/p}\chi(b_1(y)h+\dots+b_{n-1}(y)h^{n-1}+a_nh^n)\,dh\\
&=\frac{1}{p}\int_{|x|\le
1}\chi(b_1(y)px+\dots+b_{n-1}(y)p^{n-1}x^{n-1}+a_np^nx^n)\,dx,
\end{align*}
and $b_j$ is given by $(\ref{star})$.

We aim to apply Lemma~\ref{two}. When $y\neq0,$ we have
$$
|b_1(y)p|=|a_1+2a_2y+\dots+na_ny^{n-1}|/p=|ma_m|/p> p.
$$
Now if $k>j \ge 2,$ then
$$
|ma_m|\ge\left|{k-1\choose j-1}\right||ka_k|=\left|j{k\choose j}a_k\right|
$$
so that
$$
|ma_m|\ge\left|ja_j + j{j+1\choose j}a_{j+1}y+\dots+j{n\choose j}a_ny^{n-j}\right|= |jb_j(y)|.
$$
Hence if $j\ge 2,$ then
\begin{equation*}\label{you}
|jb_j(y)p^j|\le \frac{|ma_m|}{p^{j}}=\frac{|b_1(y)p|}{p^{j-1}}<|b_1(y)p|.
\end{equation*}
Thus by Lemma~\ref{two}, we have
$I_1(y)=0$ for all $y\neq0$, so that $I=I_1(0)$.
\end{proof} 
   
\noindent\textit{Proof of Lemma \ref{res}.}
We use double induction on 
$$
m=\max\{l : |la_l|\ge|ja_j| \ \textrm{ for all } \ j\neq l\},
$$ and $$r=\max_{1\le j\le
n}\log_p|ja_j|.$$
First we note trivially that
$$
|I|=\al\int_{|x|\le 1} \chi(a_1x+\dots+ a_nx^n)\,dx \ar 
\le \int_{|x|\le 1} \left|\chi(a_1x+\dots+
a_nx^n)\right|\,dx = 1.
$$
Suppose that $m=1.$ When $r\le 1,$ 
$$
\frac{p^m}{|ma_m|^{1/m}}=\frac{p}{|a_1|}\ge \frac{p}{p}= 1,
$$
so we are done. When $r>1$ we have $|a_1|>p,$ and as
$|a_1|>|ja_j|$ for all $j>1,$ we obtain the result by Lemma~\ref{two} .
Now suppose that $m>1$ and $r\le 2.$ Again we are done, as
$$
\frac{p^m}{|ma_m|^{1/m}}\ge \frac{p^2}{p^{2/2}}\ge1.
$$
So when $m=1$ or $r\le2,$ we have the result. 

Suppose we have the result when $m\le k-1$ and $r\le s-1,$ 
and suppose that $m=k$ and  $r=s.$ 
When $|y|\le 1$, we have $$\{x:|x|\le 1\}=\{x:|x-y|\le 1\}$$
so that 
\begin{align*}
|I|&=\al\int_{|x-y|\le 1} \chi(a_1x+\dots+
a_nx^n)\,dx\ar\\
&=\al\int_{|h|\le 1} \chi(a_1(h+y)+\dots+a_n(h+y)^n)d h\ar\\
&=\al\int_{|h|\le 1} \chi(b_1(y)h+\dots+ b_{n-1}(y)h^{n-1}+a_nh^n)d h\ar,
\end{align*}
for all $|y|\le 1,$ where $b_j$ is given by $(\ref{star}).$ 

As $m=k$, we have $|ka_k|>|ja_j|$ for all $j>k.$ Thus when $|y|\le1,$ 
we have $|kb_k(y)|>|jb_j(y)|$ for all $j>k,$ by Lemma~\ref{threetwo}.   
We choose $y=y_1,$ so that 
$$\max_{1\le
j<k}|jb_j(y_1)|=\min_{|y|\le 1}\max_{1\le
j<k}|jb_j(y)|.$$ Either $\max_{1\le
j<k}|jb_j(y_1)|<|kb_k(y_1)|$ or $\max_{1\le j<k}|jb_j(y_1)|\ge |kb_k(y_1)|.$ 

When $\max_{1\le
j<k}|jb_j(y_1)|<|kb_k(y_1)|,$ we have $|kb_k(y_1)|>|jb_j(y_1)|$
for all $j\neq k,$ so we can apply Lemma~\ref{four} to obtain
$$
|I|=\frac{1}{p}\left|\int_{|h|\le 1} \chi(b_1(y_1)ph+\dots+ b_{n-1}(y_1)p^{n-1}h^{n-1}+a_np^nh^n)d h\right|.
$$
Now as $\max_{1\le j\le
n}|jb_j(y_1)p^j|\le p^{s-1},$ we have $$r=\max_{1\le j \le n}\log_p|jb_j(y_1)p^j|\le s-1.$$ Since $|kb_k(y_1)|>|jb_j(y_1)|$
for all $j> k,$ we have 
$$
m=\max\{l : |lb_l(y_1)p^l|\ge|jb_j(y_1)p^j| \ \textrm{ for all } \ j\neq l\}=k_1\le k.$$ 
Hence 
$$
|I|\le   \frac{1}{p}\frac{p^{k_1}}{|k_1b_{k_1}(y_1)p^{k_1}|^{1/k_1}}=\frac{p^{k_1}}{|k_1b_{k_1}(y_1)|^{1/k_1}}\le\frac{p^{k}}{|kb_{k}(y_1)|^{1/k}},
$$
by induction. Finally 
$$|I|\le  \frac{p^k}{|ka_k|^{1/k}},$$
as $|kb_k(y_1)|=|ka_k|$ by Lemma~\ref{threetwo}.

When $\max_{1\le j<k}|jb_j(y_1)|\ge |kb_k(y_1)|$ we split the integral into $p$ pieces, so that
$$
I=\int_{|x|\le 1} \chi(a_1x+\dots+ a_nx^n)\,dx =\sum^{p-1}_{y=0}\chi(a_1y+\dots+a_ny^n)I_1(y),
$$
where 
\begin{align*}
I_1(y)&= \int_{|h|\le 1/p}\chi(b_1(y)h+\dots+b_{n-1}(y)h^{n-1}+a_nh^n)\,dh\\
&=\frac{1}{p}\int_{|x|\le
1}\chi(b_1(y)px+\dots+b_{n-1}(y)p^{n-1}x^{n-1}+a_np^nx^n)\,dx,
\end{align*}
and $b_j$ is given by $(\ref{star})$. Now when $|y|\le 1,$ we have
$$
|kb_k(y)|=|ka_k|>|lb_l(y)|
$$
for all $l>k,$ by Lemma~\ref{threetwo}. Hence
$$
\max_{1\le j<k}|jb_j(y)|\ge |kb_k(y_1)|=|ka_k|\ge |lb_l(y)|
$$
for all $l\ge k.$ Thus for $y=0,\ldots,p-1,$ there exists $k_1<k,$ where $k_1$ depends on $y$, such that $$|k_1b_{k_1}(y)|\ge|jb_j(y)|,$$ and
$$
|k_1b_{k_1}(y)p^{k_1}|>|jb_j(y)p^j|
$$
for all $j>k_1.$
Hence for $y=0,\ldots,p-1,$ we have
$$
m=\max\{l : |lb_l(y)p^l|\ge|jb_j(y)p^j| \ \textrm{ for all } \ j\neq l\}=k_1<k.$$ 
Thus
$$
|I_1(y)|\le
\frac{1}{p}\frac{p^{k_1}}{|k_1b_{k_1}(y)p^{k_1}|^{1/k_1}}=\frac{p^{k_1}}{|k_1b_{k_1}(y)|^{1/k_1}},
$$
by induction. Now as
$$
\frac{p^{k_1}}{|k_1b_{k_1}(y)|^{1/k_1}}\le\frac{p^{k-1}}{|kb_{k}(y)|^{1/k}}=\frac{p^{k-1}}{|ka_k|^{1/k}},
$$
by Lemma~\ref{threetwo}, we have
$$
|I|\le \sum_{y=0}^{p-1}|I_1(y)|\le p\frac{p^{k-1}}{|ka_k|^{1/k}}=\frac{p^{k}}{|ka_k|^{1/k}},
$$
and we are done.\hfill $\square$

\ \newline
\noindent Thanks to M. Cowling for all his invaluable help.


\begin{thebibliography}{9}


\bibitem{ar} G.I. Arhipov, A.A. Karacuba and V.N. Cubarikov, Trigonometric integrals, \textit{Math. USSR Izvestija} \textbf{15} (1980), 211--239. 
\bibitem{ca} A. Carbery, M. Christ\ and\ J. Wright, Multidimensional van der Corput and sublevel set estimates, J. Amer. Math. Soc. {\bf 12} (1999), no.~4, 981--1015. MR1683156 
\bibitem{van} J.G. van der Corput, Zahlentheoretische absch\"atzungen, \textit{Math. Ann.} \textbf{84} (1921), 53--79.
\bibitem{ge} N. Koblitz, {\it $p$-adic analysis: a short course on recent work}, Cambridge Univ. Press, Cambridge, 1980. MR0591682 
\bibitem{me} K. Rogers, Sharp van der Corput estimates and minimal divided differences, submitted.
\bibitem{st} E. M. Stein, {\it Harmonic analysis: real-variable methods, orthogonality, and oscillatory integrals}, Princeton Univ. Press, Princeton, NJ, 1993. MR1232192 
\bibitem{sa} M. H. Taibleson, {\it Fourier analysis on local fields}, Princeton Univ. Press, Princeton, N.J., 1975. MR0487295 
\bibitem{man} J. Wright, $p$-Adic van der Corput lemmas, unpublished manuscript.



\end{thebibliography}
\end{document}